\documentclass{article}
\usepackage{graphicx}
\graphicspath{{images/}}
\usepackage{amsmath, amssymb, amsthm, verbatim, comment, xcolor}
\usepackage{enumitem}
\usepackage{mathtools}
\usepackage{nccmath}
\usepackage[title]{appendix}
\newcommand{\norm}[1]{\left\lVert#1\right\rVert}
\newcommand{\ep}{_\varepsilon}
\renewcommand{\epsilon}{\varepsilon}
\newcommand{\abs}[1]{\left\lvert#1\right\rvert}
\newcommand{\pri}{^\prime}
\newcommand{\dpri}{^{\prime\prime}}
\newcommand{\too}{\longrightarrow}
\newcommand{\dd}{\,\mathrm{d}}
\usepackage{nccmath}
\usepackage{authblk}

\usepackage[margin=1in]{geometry}

\theoremstyle{plain}
\newtheorem{theorem}{Theorem}[section]
\newtheorem{proposition}[theorem]{Proposition}
\newtheorem{lemma}[theorem]{Lemma}

\theoremstyle{remark}
\newtheorem{remark}[theorem]{Remark}

\title{Regularization Independence in the One-Dimensional Schr\"odinger Equation with Delta Potential}
\author[1]{Alex Czemerinski
\thanks{alexczemerinski@ufl.edu; corresponding author}}
\author [1] {Jonathan Nunes
\thanks{jonathan.nunes@ufl.edu}}
\author[1, 2] {Sergei V. Shabanov
\thanks{shabanov@ufl.edu}}
\affil[1]{Department of Mathematics, University of Florida, Gainesville, FL 32611, USA.}
\affil[2]{Institute of Fundamental Theory, University of Florida, Gainesville, FL 32611, USA.}

\counterwithin{equation}{section}

\begin{document}

\maketitle

\begin{abstract}
We consider a class of regularizations $\{\delta\ep\}$ of the Dirac delta defined by a uniform $L^1$ bound, the
normalization $\int\delta\ep\to1$, and a truncated first-moment condition, with no restriction on sign, shape,
support, or scaling. For this class we prove that the transmission and reflection coefficients for the Hamiltonian operator
$H\ep=-\tfrac{d^2}{dx^2}-\lambda\delta\ep$ at fixed energy converge to those of the delta interaction. Under the additional assumption that each $\delta\ep$ is bounded, we show that any
family of negative eigenvalues of $H\ep$ bounded away from zero converges to $-\lambda^2/4$ when $\lambda>0$, and that no such family exists when $\lambda<0$. Both convergences are obtained with an explicit rate, expressed through two integrals of $\delta\ep$, directly computable for a given regularization. The bound-state
convergence follows from the work of Brasche, Figari, and Teta; here, our contribution is an elementary proof with an explicit rate. Their norm resolvent
convergence, however, does not by itself imply convergence of the scattering matrix, which
was previously known for families obtained by scaling a single fixed profile; we prove it
for a strictly larger class. 
\end{abstract}



\section{Introduction} \label{sec:intro}
The Dirac delta is a distribution $\delta\in\mathcal{D}'(\mathbb{R})$ defined by
\[
(\delta,\varphi):=\varphi(0),\qquad\varphi\in\mathcal{D}(\mathbb{R}),
\]
where $\mathcal{D}(\mathbb{R})=C_c^\infty(\mathbb{R})$ denotes the space of smooth compactly supported functions. The Dirac delta plays a fundamental role in the theory of generalized functions and appears frequently in exactly solvable models in quantum mechanics (see, for example,~\cite{BR}). In particular, it is used to model point interactions in the one-dimensional time-independent Schr\"odinger equation
\begin{align}\label{eq:schro}
    -\psi''(x)+V(x)\psi(x)=E\psi(x).
\end{align}
In the present work, we consider both scattering and bound-state problems. For the scattering problem, we take $E=k^2>0$ and seek a solution to \eqref{eq:schro} with a prescribed asymptotic behavior as $\abs{x}\to\infty$ corresponding to incident, transmitted, and reflected plane waves. The objective is to determine the amplitudes of the reflected and transmitted waves for an incident wave of unit amplitude. For the bound-state problem, we take $E<0$ and seek a nonzero solution $\psi\in L^2(\mathbb{R})$.

We are interested in the singular potential 
\begin{equation}\label{delta}
    V(x)=-\lambda\delta(x).
\end{equation}
Since the product $\delta\psi$ is not well defined for a general $\psi\in L^2(\mathbb{R})$, the formal expression $-\psi''-\lambda\delta\psi$ has no immediate meaning as an operator on $L^2(\mathbb{R})$. One rigorous realization of the corresponding Hamiltonian is obtained through a self-adjoint extension of $-\frac {d^2}{dx^2}$ initially defined away from the origin (see, for example,~\cite{AGHH, Gitman}). This approach produces a four-parameter family of operators, whereas only a one-parameter subfamily, characterized by specific matching conditions at the origin, is associated with the delta interaction \eqref{delta}. An alternative approach that is commonly used in physics, is to employ a \textit{regularization} of the singular potential by replacing it with a sequence of classical potentials $V\ep(x)=-\lambda\delta\ep(x)$ concentrating at the origin. The regularized potentials are used as approximations to the idealized delta interaction, with the expectation that the physical quantities should become independent of the particular choice of regularization in the limit. Most generally, a regularization of the Dirac delta potential refers to a family $\{\delta\ep\}_{\epsilon>0}$ of locally integrable functions such that
\[
\int\delta_\epsilon\varphi\to \varphi(0)\ \text{as}\ \epsilon\to0^+,\qquad\ \forall\varphi\in \mathcal{D}(\mathbb{R}).
\]
This condition alone, however, is not strong enough to guarantee the convergence of the corresponding bound-state energy or scattering amplitudes to those of the delta interaction. This raises the natural question: for what classes of regularizations $\delta\ep$ can one guarantee convergence of the spectral and scattering data to those of the delta interaction?

\subsection{The target operator}\label{ss:limitop}
Throughout the paper, $\lambda\in\mathbb{R}\setminus\{0\}$ denotes the coupling constant in the potential
$V=-\lambda\delta$, so that $\lambda>0$ is the attractive potential case.
We will compare the observables of the regularized problem against those of the standard delta interaction of strength $-\lambda$, namely the self-adjoint operator $H$ in $L^2(\mathbb{R})$ defined by
\begin{equation}\label{eq:H}
\begin{aligned}
\mathrm{dom}(H)&=\Big\{\psi\in H^2\!\left(\mathbb{R}\setminus\{0\}\right):\
\psi(0^+)=\psi(0^-),\ \psi\pri(0^+)-\psi\pri(0^-)=-\lambda\psi(0)\Big\},\\
(H\psi)(x)&=-\psi\dpri(x),\qquad x\ne 0,
\end{aligned}
\end{equation}
where $H^2$ denotes the second order Sobolev space. 
This $H$ is a member of the four-parameter family of self-adjoint extensions of $-\tfrac{d^2}{dx^2}$ restricted to
$C_c^\infty\!\left(\mathbb{R}\setminus\{0\}\right)$~\cite{Gitman}. It is often denoted by $-\Delta_{-\lambda}$~\cite{AGHH}.
The observables of $H$ easily follow from \eqref{eq:H}.
The bound-state problem admits a solution precisely when $\lambda>0$, in which case
\begin{equation}\label{eq:Estd}
E=-\frac{\lambda^2}{4},\qquad \psi(x)=\sqrt{\frac{\lambda}{2}}\;e^{-\lambda\abs{x}/2}.
\end{equation}
For scattering at fixed energy $E=k^2>0$, we seek a solution to \eqref{eq:schro} satisfying
\[
\psi(x)=
\begin{cases}
e^{ikx}+R e^{-ikx}+o(1), & x\to-\infty,\\[2pt]
T e^{ikx}+o(1), & x\to+\infty,
\end{cases}
\]
where $T$ and $R$ are called the transmission and reflection coefficients, respectively. 
Solving $-\psi\dpri=k^2\psi$ on each half-line and applying the matching conditions from \eqref{eq:H} then gives:
\begin{equation}\label{eq:TRstd}
T(k)=\frac{2k}{2k-i\lambda},\qquad R(k)=\frac{i\lambda}{2k-i\lambda}.
\end{equation}
The essential spectrum is $\sigma_{\mathrm{ess}}(H)=[0,\infty)$, and the point spectrum contains a single value as given in \eqref{eq:Estd} ~\cite[Thm.~I.3.1.4]{AGHH}.
Note that $H$ enters this paper only through the target values \eqref{eq:Estd} and \eqref{eq:TRstd}. We do not study operator convergence of $-\frac{d^2}{dx^2}-\lambda\delta\ep$ to $H$.

\subsection{Prior results}\label{sec:priorResults}
There is a substantial literature on approximations of point interactions in Schr\"odinger operators~\cite{AGHH,BFT, GH2010, GH2013, Friedman1972,AHK1981,AGHK1984,
AlbeverioNizhnik2000,Golovaty2012,
DiasJorgePrata2016}. Among the results most relevant to the present work are those of Albeverio et al.~\cite{AGHH} and Brasche, Figari, and Teta~\cite{BFT}.

Albeverio et al.\ in~\cite[Ch.~I.3]{AGHH} treat families obtained by scaling a single fixed function. For real-valued $V\in L^1(\mathbb{R})$ and $\gamma(\varepsilon)$ real-analytic near the origin with $\gamma(0)=0$, the operators $-\frac{d^2}{dx^2}+\gamma(\varepsilon)\varepsilon^{-2}V(x/\epsilon)$ converge in the norm resolvent sense to the delta interaction of strength $\alpha=\gamma\pri(0)\int V$~\cite[Thm.~I.3.2.3]{AGHH}, and the corresponding scattering matrix converges to that of the limiting operator~\cite[Thm.~I.3.4.1]{AGHH}.

Brasche, Figari, and Teta~\cite{BFT} consider operators $-\Delta+\mu$, where
$\mu$ is a finite signed Radon measure on $\mathbb{R}$, defined via the
associated quadratic form. Their Theorem~3 shows that weak
convergence of such measures, meaning that $\int f\,d\mu_n\to\int f\,d\mu$ for every
$f\in C_b(\mathbb{R})$, the space of bounded continuous functions, implies the
norm resolvent convergence of $-\Delta+\mu_n$ to $-\Delta+\mu$. The proof also
shows that the corresponding quadratic forms are bounded below uniformly for all sufficiently large $n$.
Norm resolvent convergence implies $\sigma(-\Delta+\mu)=\lim_{n\to\infty}
\sigma(-\Delta+\mu_n)$~\cite[Thm.~6.38]{Teschl}, so every finite limit point of
eigenvalues of $-\Delta+\mu_n$ lies in $\sigma(-\Delta+\mu)$. This fact, together with the uniform lower bound to discard eigenvalues tending to $-\infty$, settles the convergence of any family of negative eigenvalues bounded away from zero to the ground-state energy of the delta interaction. This is true for any regularization whose associated measures converge weakly to $-\lambda$ times the Dirac measure, which includes the families considered here. No rate of convergence is given in~\cite{BFT}, and scattering data are not
treated there.

Related works also include that of Golovaty and Hryniv~\cite{GH2010,GH2013}, who study regularizations
of the $\delta'$-potential. They prove norm resolvent convergence for single-shape scaled families
$\varepsilon^{-2}V(x/\varepsilon)$, first for compactly supported $V$~\cite{GH2010} and then for $V$
in the Faddeev--Marchenko class~\cite{GH2013}, with a limiting operator that depends on whether $V$
possesses a half-bound state.

\subsection{Main results}\label{ss:results}
Let $\{\delta\ep\}_{\epsilon>0}$ be a family of real-valued integrable functions on $\mathbb{R}$. Put
\begin{equation}\label{eq:mp}
m\ep:=\int\delta\ep(x)\dd x,
\qquad
p\ep:=\int\min\{1,\abs{x}\}\abs{\delta\ep(x)}\dd x,
\qquad
r\ep:=p\ep+\abs{m\ep-1}.
\end{equation}
We refer to the numbers $m\ep$ and $p\ep$ as the total mass and the truncated moment, respectively. Suppose that
\[
\int\abs{\delta\ep}\dd x\le M\ \text{ for some }M\text{ independent of }\epsilon,
\qquad m\ep\to1,
\quad p\ep\to0, \quad \text{as}\quad \epsilon\to 0^+ .
\]
Thus $\delta\ep$ is required to have a uniformly bounded $L^1$ norm, to have the total mass tending to one, and to
concentrate at the origin in the sense measured by the moment $p\ep$. No sign condition, no fixed
shape, and no compact support is imposed.

The following results are established in this paper. First, the scattering coefficients $T\ep,R\ep$ at fixed $k>0$ converge to \eqref{eq:TRstd}
(Theorem~\ref{thm:scattering}). Second, if in addition each $\delta\ep\in L^\infty(\mathbb{R})$, a family of negative eigenvalues $E\ep$ bounded away from $0$ can exist only when $\lambda>0$, in which case $E\ep\to-\lambda^2/4$ (Theorem~\ref{thm:bound}). Third, both convergences carry the explicit rate $O(r\ep)$ with computable constants.

The two proofs follow a similar structure. In both the scattering and bound-state problems, concentration of $\delta\ep$ allows the quantities appearing in the corresponding integral equation to be replaced by their values at the origin, up to an error controlled by $p\ep$. The resulting scalar relation is then used to determine the limiting scattering coefficients or bound-state energy.

\subsection{Relation to prior work}\label{ss:contribution}
A precise description of the families of regularizations considered in this paper is given in Section~\ref{ss:class} (see hypotheses \ref{A1}--\ref{A4}), and all such regularizations satisfy the hypotheses of~\cite[Thm.~3]{BFT}. Indeed, $\dd\mu\ep:=-\lambda\delta\ep\dd x$
defines a finite signed Radon measure by \ref{A1}, and \ref{A1}--\ref{A3} imply that $\mu\ep$ converges weakly to
$-\lambda$ times the Dirac measure at the origin (see Appendix~\ref{app:weak}). Consequently, the associated operators converge to the delta interaction in the norm resolvent sense. Since norm
resolvent convergence~\cite[Thm.~6.38]{Teschl}, together with the uniform lower
bound established in the proof of~\cite[Thm.~3]{BFT}, implies convergence of
eigenvalues bounded away from zero, the qualitative part of our
Theorem~\ref{thm:bound} follows.

What the present paper adds can be summarized as follows. The principal new
result is the extension of the convergence of scattering data to the broader class
\ref{A1}--\ref{A3}. In addition, we provide an elementary proof and obtain a quantitative rate for the convergence of the scattering coefficients and bound-state energy.

\begin{enumerate}[label=\textup{(\roman*)},leftmargin=2.4em]

\item \emph{Scattering data convergence for a broader class.}
Albeverio et al.\ establish convergence of the transmission and reflection
coefficients for families obtained by scaling a single fixed profile
$V\in L^1$~\cite[Thm.~I.3.4.1]{AGHH}. Every such family with a nontrivial limit, divided by its strength $\alpha=\gamma\pri(0)\int V\neq0$, satisfies
\ref{A1}--\ref{A3}, so our Theorem~\ref{thm:scattering} with $\lambda=-\alpha$ recovers their convergence result, while the converse fails because
\ref{A1}--\ref{A3} allow the shape to vary with $\epsilon$.
Thus Theorem~\ref{thm:scattering} extends this result to a strictly broader
class of regularizations. No scattering theory is developed in~\cite{BFT},
and convergence of the scattering coefficients is not a consequence of
norm resolvent convergence alone.

\item \emph{Quantitative convergence rates.}
Theorem~3 of~\cite{BFT} is qualitative and gives no rate of convergence.
Our proof yields
\[
\abs{T_\epsilon-T}
\leq A_k\abs{m\ep-1}+B_kp\ep,
\]
where the constants $A_k$ and $B_k$ are explicit and depend only on $k, \lambda, $ and $M$. We obtain similar expressions for the
reflection coefficients (Theorem~\ref{thm:scattering}) and bound-state energies (Theorem~\ref{thm:bound}). Since $p_\epsilon$ and $m_\epsilon$
are defined directly in terms of $\delta_\epsilon$, the rates of convergence can be
computed explicitly for a given regularization.

\item \emph{An elementary proof of the bound-state convergence.}
Although the bound-state convergence follows from~\cite[Thm.~3]{BFT},
our proof is direct and uses only the integral-equation representation and
elementary inequalities.

\end{enumerate}

\section{Preliminaries}
\subsection{The class of regularizations}\label{ss:class}

The family $\{\delta\ep\}_{\epsilon>0}$ considered here consists of real-valued measurable functions satisfying \ref{A1}--\ref{A3}.
\begin{enumerate}[label=\textup{(A\arabic*)},leftmargin=3.2em]
\item\label{A1} $\medint\int\abs{\delta\ep(x)}\dd x\leq M$ for a constant $M$ independent of $\epsilon$;  
\item\label{A2} $m\ep\to 1$ as $\epsilon\to 0^+$;
\item\label{A3} $p\ep \to0$ as $\epsilon\to0^+$,
\end{enumerate}
where $m\ep$ and $p\ep$ are as defined in \eqref{eq:mp}. In Section~\ref{sec:bound} we assume, in addition, that
\begin{enumerate}[label=\textup{(A\arabic*)},leftmargin=3.2em,start=4]
\item\label{A4} $\delta\ep\in L^\infty(\mathbb{R})$ for every $\epsilon>0$.
\end{enumerate}

\subsection{Elementary lemma}

\begin{lemma}[Elementary estimates]\label{lem:elem}
Let $s,t,\theta\in\mathbb{R}$, $v>0$, $f\in H^1(\mathbb{R})$, the first order Sobolev space. Then
\begin{enumerate}[label=\textup{(\alph*)},leftmargin=2.4em]
\item $\min\{1,\abs{s-t}\}\le\min\{1,\abs{s}\}+\min\{1,\abs{t}\}$;
\item $\abs{e^{i\theta t}-1}\le 2(1+\abs{\theta})\min\{1,\abs{t}\}$;
\item $\abs{e^{-v\abs{t}}-1}\le(1+v)\min\{1,\abs{t}\}$;
\item $\norm{f}_\infty^2\le2\norm{f\pri}_2\norm{f}_2$.
\end{enumerate}
\end{lemma}
\begin{proof}
(a) Let $\phi(u)=\min{\{1, u\}}$. If $\abs{s}\geq1$ or $\abs{t}\geq1$, then $\phi(\abs{s})+\phi(\abs{t})\geq1\geq\phi(\abs{s-t})$. Otherwise, $\abs{s},\abs{t}<1$ so $\phi(\abs{s})+\phi(\abs{t})=\abs{s}+\abs{t}\geq\abs{s-t}\geq\phi(\abs{s-t})$.
\\(b) $\abs{e^{i\theta t}-1}=2\abs{\sin(\theta t/2)}\le\min\{2,\abs{\theta}\abs{t}\}\leq2(1+\abs{\theta})\min{\{1,\abs{t}\}}$.
\\(c) $\abs{e^{-v\abs{t}}-1}\le\min\{1,v\abs{t}\}\leq(1+v)\min{\{1,\abs{t}\}}$.
\\(d) $\abs{f(x)}^2=\int_{-\infty}^x\left(\abs{f}^2\right)\pri\le2\int\abs{f}\abs{f\pri}\le2\norm{f}_2\norm{f\pri}_2$.
\end{proof}

\section{Regularization Independence for the Scattering Problem}\label{sec:scattering}

Fix $k>0$ and consider the scattering problem with a regularized delta potential:
\begin{equation}\label{eq:scat}
-\psi\ep\dpri(x)-\lambda\delta\ep(x)\psi\ep(x)=k^2\psi\ep(x),
\qquad
\psi\ep(x)=
\begin{cases}
e^{ikx}+R\ep e^{-ikx}+o(1), & x\to-\infty,\\[2pt]
T\ep e^{ikx}+o(1), & x\to+\infty,
\end{cases}
\end{equation}
where the coefficients $T\ep,R\ep\in\mathbb{C}$ are determined by the asymptotic behavior of the solution. Here $o(1)$ denotes terms that tend to zero as $x\to\pm\infty$ for fixed $\epsilon>0$.

Since $\delta\ep$ is merely integrable, the differential equation in \eqref{eq:scat} is understood in the sense of distributions. For the product $\delta\ep\psi\ep$ to define a distribution, we require
\[
\psi\ep\in L^1_{\mathrm{loc}}(\mathbb{R})\qquad\text{and}\qquad\delta\ep\psi\ep\in L^1_{\mathrm{loc}}(\mathbb{R}),
\]
and we impose no further conditions on $\psi\ep$.
\begin{remark}\label{rem:existence}
For every real-valued $\delta\ep\in L^1(\mathbb{R})$ and every $k>0$, problem
\eqref{eq:scat} has a solution in the above sense. Indeed, it is shown in
\cite[Cor.~2.5, Lemma~2.6, Eq.~(2.15)]{HM2021} that the equation in
\eqref{eq:scat} has a continuous solution $\phi$ satisfying
\[
\phi(x)=e^{ikx}+o(1) \text{ as } x\to+\infty
\quad\text{and}\quad
\phi(x)=ae^{ikx}+be^{-ikx}+o(1) \text{ as } x\to-\infty,
\]
for some $a,b\in\mathbb{C}$ with $|a|^2-|b|^2=1$. In particular $a\neq0$, so
$\psi\ep:=\phi/a$ solves \eqref{eq:scat} with $T\ep=1/a$ and $R\ep=b/a$.
\end{remark}
\begin{theorem}[Regularization independence of the scattering coefficients]\label{thm:scattering}
Let $\{\delta\ep\}$ satisfy \ref{A1} and $k>0$. For each $\epsilon>0$, let
$\psi\ep$ be a solution to \eqref{eq:scat}, with coefficients $T\ep$ and $R\ep$, and
\[
T=\frac{2k}{2k-i\lambda},\qquad R=\frac{i\lambda}{2k-i\lambda}
\]
be the coefficients \eqref{eq:TRstd} of the delta interaction. Then,
\[
\abs{T\ep-T}\le A_k\abs{m\ep-1}+B_kp\ep,
\qquad
\abs{R\ep-R}\le A_k\abs{m\ep-1}+B_kp\ep,
\]
where
\[
C_k:=\sqrt{2}e^{M\abs{\lambda}/2k},\qquad
A_k:=\frac{\abs{\lambda}}{2k}\left(1+\frac{\abs{\lambda}MC_k}{2k}\right),\qquad
B_k:=\frac{\abs{\lambda}}{k}(1+k)\left(1+C_k\left(1+\frac{\abs{\lambda}M}{k}\right)\right).
\]
In particular, if in addition \ref{A2} and \ref{A3} hold, then $T\ep\to T$ and $R\ep\to R$ as $\epsilon\to0^+$, with rate $O(r\ep)$.
\end{theorem}

\subsection{Proof of Theorem~\ref{thm:scattering}}
We begin with the regularity analysis of the solutions to \eqref{eq:scat}. The equation gives $\psi\ep\dpri=-(\lambda\delta\ep+k^2)\psi\ep\in L^1_{\mathrm{loc}}(\mathbb{R})$, in the sense of distributions. Hence, we can pick a representative such that $\psi\ep\pri\in AC_{\mathrm{loc}}$, the space of locally absolutely continuous functions. In particular, this gives $\psi\ep\in C^1(\mathbb{R})$. Continuity implies that $\psi\ep$ is locally bounded; together with the asymptotics in \eqref{eq:scat}, this gives $\psi\ep\in L^\infty(\mathbb{R})$. Since \ref{A1} gives $\delta\ep\in L^1$, we have $\delta\ep\psi\ep\in L^1(\mathbb{R})$.

Next, we show that every solution of \eqref{eq:scat} satisfies the
Lippmann--Schwinger equation, which expresses $\psi\ep$ as the incoming wave
plus the convolution of $\lambda\delta\ep\psi\ep$ with the outgoing Green's
function $\frac{i}{2k}e^{ik|x|}$ of $-\frac{d^2}{dx^2}-k^2$.

\begin{lemma}[Lippmann--Schwinger equation]\label{lem:LS}
Let $\delta\ep$ satisfy \ref{A1} and let $\psi\ep$ solve \eqref{eq:scat}. Then, for every $x\in\mathbb{R}$,
\begin{equation}\label{eq:Lippmann}
\psi\ep(x)=e^{ikx}+\frac{i\lambda}{2k}\int e^{ik|x-y|}
\delta\ep(y)\psi\ep(y)\,\mathrm{d}y .
\end{equation}
\end{lemma}
The proof is standard. For completeness, we give it in Appendix~\ref{app:LS}, using only that $\delta\ep\in L^1(\mathbb{R})$. The proof also yields the asymptotic behavior of $\psi\ep'$:
\begin{equation}\label{eq:psiprime}
\psi\ep'(x)=
\begin{cases}
ik\bigl(e^{ikx}-R\ep e^{-ikx}\bigr)+o(1), & x\to-\infty,\\[2pt]
ik\,T\ep e^{ikx}+o(1), & x\to+\infty,
\end{cases}
\end{equation}
which we use in the proof of Proposition~\ref{prop:S1}
\newline
Define 
\[c\ep:=\int\delta\ep(y)\psi\ep(y)\dd y.\] The conclusion of Theorem~\ref{thm:scattering} follows from \eqref{eq:Lippmann} and the following three propositions, proved in the subsections below.  

\begin{proposition}[Uniform bound]\label{prop:S1}
Assume \ref{A1} and let $\psi\ep$ solve \eqref{eq:scat}. Then $\norm{\psi\ep}_\infty\le C_k$ for every $\epsilon>0$, with $C_k$ as in Theorem~\ref{thm:scattering}.
\end{proposition}

\begin{proposition}[Coefficients in terms of $c\ep$]\label{prop:S2}
Assume \ref{A1} and let $\psi\ep$ solve \eqref{eq:scat}. Then, for every $\epsilon>0$,
\[
\abs{T\ep-1-\frac{i\lambda}{2k}c\ep}\le\frac{\abs{\lambda}}{k}(1+k)C_kp\ep,
\qquad
\abs{R\ep-\frac{i\lambda}{2k}c\ep}\le\frac{\abs{\lambda}}{k}(1+k)C_kp\ep.
\]
\end{proposition}

\begin{proposition}[Scalar equation]\label{prop:S3}
Assume \ref{A1} and let $\psi\ep$ solve \eqref{eq:scat}. Then, for every $\varepsilon>0$,
\[
\abs{\frac{i\lambda}{2k}c\ep-\frac{i\lambda}{2k-i\lambda}}
\le A_k|m\ep-1|+\frac{|\lambda|}{k}(1+k)\Bigl(1+\frac{|\lambda|MC_k}{k}\Bigr)p\ep ,
\]
with $A_k$ as in Theorem~\ref{thm:scattering}.
\end{proposition}
We now prove Theorem~\ref{thm:scattering}, assuming Propositions~\ref{prop:S1}--\ref{prop:S3}.
Since $T=\medmath {\frac{2k}{2k-i\lambda}=1+\frac{i\lambda}{2k-i\lambda}}$ and $R=\medmath{\frac{i\lambda}{2k-i\lambda}}$, the triangle inequality gives
\[
|T\ep-T|\le\abs{T\ep-1-\frac{i\lambda}{2k}c\ep}+\abs{\frac{i\lambda}{2k}c\ep-\frac{i\lambda}{2k-i\lambda}},
\qquad
\abs{R\ep-R}\le\abs{R\ep-\frac{i\lambda}{2k}c\ep}+\abs{\frac{i\lambda}{2k}c\ep-\frac{i\lambda}{2k-i\lambda}},
\]
and Propositions~\ref{prop:S2} and~\ref{prop:S3} give the bounds of Theorem~\ref{thm:scattering}, since the coefficients of $p\ep$ add up to $B_k$.
Finally, under \ref{A2} and \ref{A3}, the convergence follows. \hfill \qedsymbol

\subsection{Proof of Proposition~\ref{prop:S1}}

Define
\[
\mathcal{E}(x):=\abs{\psi\ep(x)}^2+\frac{1}{k^2}\abs{\psi\pri\ep(x)}^2 .
\]
Since $\psi\ep,\psi\ep\pri\in AC_{\mathrm{loc}}$, so is $\mathcal{E}$. Differentiating and substituting $\psi\ep\dpri=\left(-\lambda\delta\ep-k^2\right)\psi\ep$, two terms cancel and
\[
\mathcal{E}\pri(x)=2\Re\left(\psi\ep\overline{\psi\ep\pri}\right)
+\frac{2}{k^2}\Re\left(\overline{\psi\ep\pri}\left(-\lambda\delta\ep-k^2\right)\psi\ep\right)
=\frac{-2\lambda\delta\ep}{k^2}\Re\left(\overline{\psi\ep\pri}\psi\ep\right).
\]
Hence, using $2ab\le a^2+b^2$,
\[
\abs{\mathcal{E}\pri(x)}\leq \frac{2\abs{\lambda}\abs{\delta\ep}}{k^2}\abs{\psi\ep\pri}\abs{\psi\ep}
=\frac{\abs{\lambda}\abs{\delta\ep}}{k}\left[2\abs{\psi\ep}\cdot\frac{1}{k}\abs{\psi\ep\pri}\right]
\leq\frac{\abs{\lambda}\abs{\delta\ep}}{k}\mathcal{E}(x).
\]
Note that $\mathcal{E}>0$ everywhere. Indeed, if $\mathcal{E}(x_0)=0$, then
$\psi\ep(x_0)=\psi\ep\pri(x_0)=0$, and uniqueness for the initial value problem with locally integrable
coefficients forces $\psi\ep\equiv0$, which contradicts the asymptotic condition in
\eqref{eq:scat}. Therefore, for $b>x$,
\[
\abs{\ln (\mathcal{E}(b))-\ln (\mathcal{E}(x))}=\abs{\int_x^b\frac{\mathcal{E}\pri(y)}{\mathcal{E}(y)}\dd y}
\leq\int_x^b\frac{\abs{\lambda}\abs{\delta\ep(y)}}{k}\dd y\leq\frac{M\abs{\lambda}}{k},
\]
by \ref{A1}. Hence $\mathcal{E}(x)\le \mathcal{E}(b)e^{M\abs{\lambda}/k}$. By \eqref{eq:scat} and \eqref{eq:psiprime}, $\mathcal{E}(b)\to2\abs{T\ep}^2$ as $b\to\infty$. Moreover, $\abs{T\ep}\le1$ by the flux conservation: since $\delta\ep$ is real-valued, $\overline{\psi\ep}$ also solves the differential equation in \eqref{eq:scat}, so the Wronskian $W:=\overline{\psi\ep}\psi\ep\pri-\overline{\psi\ep\pri}\psi\ep$ is constant (it is locally absolutely continuous, and $W\pri=\overline{\psi\ep}\psi\ep\dpri-\overline{\psi\ep\dpri}\psi\ep=0$ a.e.), and by \eqref{eq:scat} and \eqref{eq:psiprime},
\[
W(x)=
\begin{cases}
2ik\big(1-\abs{R\ep}^2\big)+o(1), & x\to-\infty,\\[2pt]
2ik\abs{T\ep}^2+o(1), & x\to+\infty,
\end{cases}
\]
so that $\abs{T\ep}^2+\abs{R\ep}^2=1$. Letting $b\to\infty$ therefore gives $\mathcal{E}(x)\le2e^{M\abs{\lambda}/k}$. Since $\abs{\psi\ep}^2\le \mathcal{E}$, the claim follows. \hfill \qedsymbol

\subsection{Proof of Proposition~\ref{prop:S2}}
By Proposition~\ref{prop:S1} and \ref{A1} the integral $c\ep$ is finite, with $\abs{c\ep}\le MC_k$.
Using an explicit form of the absolute value $\abs{x-y}$, \eqref{eq:Lippmann} assumes the form
\begin{equation}\label{eq:splitting}
\psi\ep(x)=e^{ikx}+\frac{i\lambda}{2k}\left[e^{ikx}\int_{-\infty}^x e^{-iky}\delta\ep(y)\psi\ep(y)\dd y
+e^{-ikx}\int_x^\infty e^{iky}\delta\ep(y)\psi\ep(y)\dd y\right].
\end{equation}
Since $\delta\ep\psi\ep\in L^1$, the tail $\medmath\int_x^\infty\abs{\delta\ep}\abs{\psi\ep}$ tends to $0$ as
$x\to+\infty$. Multiplying \eqref{eq:splitting} by $e^{-ikx}$ and letting $x\to+\infty$, the left side tends to $T\ep$ by
\eqref{eq:scat} while the last term tends to $0$, so that
\[
T\ep=1+\frac{i\lambda}{2k}\int e^{-iky}\delta\ep(y)\psi\ep(y)\dd y
=1+\frac{i\lambda}{2k}c\ep+\frac{i\lambda}{2k}\int\left(e^{-iky}-1\right)\delta\ep(y)\psi\ep(y)\dd y.
\]
Similarly, subtracting $e^{ikx}$ from both sides of \eqref{eq:splitting}, then multiplying the result by $e^{ikx}$ and letting
$x\to-\infty$, we obtain the following expression for $R\ep$,
\[
R\ep=\frac{i\lambda}{2k}\int e^{iky}\delta\ep(y)\psi\ep(y)\dd y=\frac{i\lambda}{2k} c\ep+\frac{i\lambda}{2k}\int (e^{iky}-1)\delta\ep(y)\psi\ep(y)\dd y.
\]
By Lemma~\ref{lem:elem}(b) and Proposition~\ref{prop:S1},
\[
\abs{\frac{i\lambda}{2k}\int\left(e^{\pm iky}-1\right)\delta\ep(y)\psi\ep(y)\dd y}\leq\frac{\abs{\lambda}}{k}(1+k)C_k\int\abs{\delta\ep(y)}\min{\{1,\abs{y}\}}\dd y=\frac{\abs{\lambda}}{k}(1+k)C_kp\ep,
\]
which gives the claimed bounds for $T\ep$ and $R\ep$. \hfill \qedsymbol

\subsection{Proof of Proposition~\ref{prop:S3}}
Multiplication of  \eqref{eq:Lippmann} by $\delta\ep(x)$ and subsequent integration yield
\begin{equation}\label{eq:ceps-scat}
    c\ep=I\ep+\frac{i\lambda}{2k}J\ep,
\end{equation}
where
\[
I\ep:=\int\delta\ep(x)e^{ikx}\dd x,
\qquad
J\ep:=\iint e^{ik\abs{x-y}}\delta\ep(x)\delta\ep(y)\psi\ep(y)\dd y\dd x ,
\]
both of which are finite by \ref{A1} and Proposition~\ref{prop:S1}.

By Lemma~\ref{lem:elem}(b), $\abs{I\ep-m\ep}\le2(1+k)p\ep$, so
\[
\abs{I\ep-1}\le\abs{m\ep-1}+2(1+k)p\ep .
\]
Next,
\[
J\ep=m\ep c\ep+\iint\left(e^{ik\abs{x-y}}-1\right)\delta\ep(x)\delta\ep(y)\psi\ep(y)\dd y\dd x,
\]
and applying Lemma~\ref{lem:elem}(b), then Lemma~\ref{lem:elem}(a) to split $\min\{1,\abs{x-y}\}$, then \ref{A1},
\[
\abs{J\ep-m\ep c\ep}\le2(1+k)C_k\iint\min\{1,\abs{x-y}\}\abs{\delta\ep(x)}\abs{\delta\ep(y)}\dd y\dd x
\le4(1+k)C_kMp\ep .
\]
By \eqref{eq:ceps-scat}, and the identity $J\ep=c\ep+(m\ep-1)c\ep+(J\ep-m\ep c\ep)$, we obtain
\[
c\ep=1+\frac{i\lambda}{2k}c\ep+e\ep,\qquad
e\ep:=(I\ep-1)+\frac{i\lambda}{2k}\bigl[(J\ep-m\ep c\ep)+(m\ep-1)c\ep\bigr].
\]
Using the two estimates above, together with $|c\ep|\le MC_k$, we infer that
\[
\abs{e\ep}\leq\Bigl(1+\frac{|\lambda|MC_k}{2k}\Bigr)|m\ep-1|+2(1+k)\Bigl(1+\frac{|\lambda|MC_k}{k}\Bigr)p\ep .
\]
Multiplying the equation for $c\ep$ by $\frac{i\lambda}{2k}$ and solving it for $\frac{i\lambda}{2k}c\ep$ gives
$\frac{i\lambda}{2k}c\ep=\frac{i\lambda}{2k-i\lambda}(1+e\ep)$. Since
$\abs{\frac{i\lambda}{2k-i\lambda}}\le\frac{|\lambda|}{2k}$, the claimed bound follows. \hfill \qedsymbol

\section{Regularization Independence for the Bound-State Problem}\label{sec:bound}

Consider the Schr\"odinger equation with a regularized delta potential at energy $E\ep=-v\ep^2$, $v\ep>0$:
\begin{align}\label{eq:bound}
    -\psi_\epsilon''(x)-\lambda\delta\ep(x)\psi\ep(x)=-v\ep^2\psi\ep(x),
    \quad \psi_\epsilon\in L^2,\quad\norm{\psi\ep}_2=1.
\end{align}
\begin{theorem}[Regularization independence of the bound state]\label{thm:bound}
Let $\{\delta\ep\}$ satisfy \ref{A1}--\ref{A4} and $\lambda>0$, and suppose that for all
sufficiently small $\varepsilon>0$ the pair $(v\ep,\psi\ep)$ solves \eqref{eq:bound}, with
$\liminf_{\varepsilon\to0^+}v\ep>0$. Then, for all sufficiently small $\varepsilon$,
\[
\abs{v\ep-\frac{\lambda}{2}}\le\frac{\lambda}{2}\,|m\ep-1|+B\,p\ep,
\qquad\text{where}\qquad B:=32\lambda M^3(\lambda+12+8\lambda M).
\]
In particular, as $\varepsilon\to0^+$,
\[
v\ep\to\frac{\lambda}{2},\qquad\text{equivalently}\qquad E\ep\to-\frac{\lambda^2}{4},
\]
the bound-state energy \eqref{eq:Estd} of the delta interaction. If instead $\lambda<0$, there
is no family $\{(v\ep,\psi\ep)\}_{\varepsilon>0}$ of solutions of \eqref{eq:bound} with
$\liminf_{\varepsilon\to0^+}v\ep>0$.
\end{theorem}
\subsection{Proof of Theorem~\ref{thm:bound}}
Let $G_{v\ep}$ be the Green's function of $-\frac{d^2}{dx^2}+v\ep^2$ that decays as $\abs{x}\to\infty$, namely
\begin{align}
    G_{v\ep}(x)&=\frac{1}{2v\ep}e^{-v\ep\lvert x\rvert}. \label{greens}
\end{align}
\begin{lemma}[Convolution equation]\label{lem:conv}
Let $\delta\ep$ satisfy \ref{A4} and let $(v\ep,\psi\ep)$ solve \eqref{eq:bound}. Then $\psi\ep\in H^2(\mathbb{R})$ and for every $x\in\mathbb{R}$,
\begin{equation}\label{conv}
    \psi\ep(x)=\lambda\big[G_{v\ep}*(\delta\ep\psi\ep)\big](x)=\frac\lambda{2v\ep}\int e^{-v\ep\abs{x-y}}\delta\ep(y)\psi\ep(y)\dd y .
\end{equation}
\end{lemma}
\begin{proof}
By \ref{A4}, $\delta\ep\psi\ep\in L^2(\mathbb{R})$, so \eqref{eq:bound} gives $\psi\ep\dpri=v\ep^2\psi\ep-\lambda\delta\ep\psi\ep\in L^2(\mathbb{R})$. Since also $\psi\ep\in L^2(\mathbb{R})$, we have $\psi\ep\in H^2(\mathbb{R})$. The operator $H\ep=-\frac{d^2}{dx^2}+v\ep^2\colon H^2(\mathbb{R})\to L^2(\mathbb{R})$ is invertible, and its inverse is convolution with $G_{v\ep}$ (see, for example,~\cite[\S7.5]{Teschl}). Applying it to $\big(-\frac{d^2}{dx^2}+v\ep^2\big)\psi\ep=\lambda\delta\ep\psi\ep$ gives \eqref{conv} in $L^2(\mathbb{R})$. Both sides of \eqref{conv} are continuous, the left one by the Sobolev embedding theorem $H^2(\mathbb{R})\hookrightarrow C^1(\mathbb{R})$ and the right one since $G_{v\ep},\delta\ep\psi\ep\in L^2(\mathbb{R})$, so \eqref{conv} holds for every $x\in\mathbb{R}$.
\end{proof}

Before proving Theorem~\ref{thm:bound}, we state three technical facts, proved in the subsections below.
Fix $v_{\min}$ with $0<v_{\min}<\liminf_{\epsilon\to0^+}v\ep$, so that for all sufficiently small $\epsilon$,
\[
v\ep\geq v_{\min}.
\]

\begin{proposition}[Properties of the solution]\label{prop:B1}
Assume \ref{A1}--\ref{A4} and suppose $(v\ep,\psi\ep)$ solves \eqref{eq:bound}. Then, for all sufficiently small $\epsilon$, the following bounds hold:
\begin{enumerate}[label=\textup{(\alph*)},leftmargin=2.4em]
\item $\displaystyle\frac{v_{\min}}{\sqrt{\abs{\lambda}M}}\leq\norm{\psi\ep}_\infty\leq 2\sqrt{\abs{\lambda}M}$,
\item $\norm{\psi\pri\ep}_\infty\leq (\abs{\lambda}M)^{3/2}$,
\item $v\ep\leq2\abs{\lambda}M$.
\end{enumerate}
\end{proposition}

\begin{proposition}[Error terms vanish]\label{prop:B2}
Assume \ref{A1}--\ref{A4} and suppose $(v\ep,\psi\ep)$ solves \eqref{eq:bound}. Then, for all sufficiently small $\epsilon$,
\begin{enumerate}[label=\textup{(\alph*)},leftmargin=2.4em]
\item $\displaystyle\frac{\abs{\lambda}}{2v\ep}\int \abs{e^{-v\ep\abs{y}}-1}\abs{\delta\ep(y)}\abs{\psi\ep(y)}\dd y\leq C_1p\ep$,
\item $\displaystyle\frac{\abs{\lambda}}{2v\ep}\int \abs{\delta\ep(y)}\abs{\psi\ep(y)-\psi\ep(0)}\dd y\leq C_2p\ep$,
\end{enumerate}
where
\[
C_1:=\abs{\lambda}\sqrt{\abs{\lambda}M}\left(1+\frac{1}{v_{\min}}\right),\qquad
C_2:=\frac{2}{v_{\min}}\abs{\lambda}\sqrt{\abs{\lambda}M}(1+\abs{\lambda}M).
\]
\end{proposition}

\begin{proposition}[Lower bound on $\psi\ep(0)$]\label{prop:B3}
Assume \ref{A1}--\ref{A4} and suppose $(v\ep,\psi\ep)$ solves \eqref{eq:bound}. Then, for all sufficiently small $\epsilon$,

\[
\abs{\psi\ep(0)}\geq \frac{{v_{\min}}^2}{(\abs{\lambda}M)^{3/2}}. 
\]
\end{proposition}
We now prove Theorem~\ref{thm:bound}, assuming Propositions~\ref{prop:B1}--\ref{prop:B3}. By Lemma~\ref{lem:conv}, $\psi\ep$ is continuous. Evaluating \eqref{conv} at $x=0$, we obtain
\begin{align*}
    \psi\ep(0)&=\frac\lambda{2v\ep}\int e^{-v\ep\abs{y}}\delta\ep(y)\psi\ep(y)\dd y,
    \\&=\frac\lambda{2v\ep}\int (e^{-v\ep\abs{y}}-1)\delta\ep(y)\psi\ep(y)\dd y+\frac\lambda{2v\ep}\int \delta\ep(y)\psi\ep(y)\dd y,
    \\&=\frac\lambda{2v\ep}\left[\int (e^{-v\ep\abs{y}}-1)\delta\ep(y)\psi\ep(y)\dd y+\int \delta\ep(y)(\psi\ep(y)-\psi\ep(0))\dd y+\int \delta\ep(y)\psi\ep(0)\dd y\right].
\end{align*}
By Proposition~\ref{prop:B2}, the absolute values of the first two terms in the last equality are bounded by $(C_1+C_2)p\ep$. Thus,
\[
\psi\ep(0)=\frac\lambda{2v\ep}m\ep\psi\ep(0)+e\ep,\qquad\abs{e\ep}\le(C_1+C_2)p\ep.
\]
Since $\psi\ep(0)\neq0$ by Proposition~\ref{prop:B3}, dividing this equation by $\psi\ep(0)$ and multiplying it by $v\ep$ gives $v\ep-\frac\lambda2m\ep=v\ep e\ep/\psi\ep(0)$. Hence, by Propositions~\ref{prop:B1}(c) and \ref{prop:B3},
\begin{align*}
\abs{v\ep-\frac\lambda2m\ep}\le2\abs{\lambda}M(C_1+C_2)\frac{(\abs{\lambda}M)^{3/2}}{v_{\min}^2}p\ep=B_{v_{\min}}p\ep,
\end{align*}
where $B_{v_{\min}}:=2\lambda^4M^3(v_{\min}+3+2\abs{\lambda}M)/v_{\min}^3$. Hence, for all sufficiently small $\epsilon$,
\[
\abs{v\ep-\frac\lambda2}\le\abs{v\ep-\frac\lambda2m\ep}+\frac{\abs{\lambda}}{2}\abs{m\ep-1}\le\frac{\abs{\lambda}}{2}\abs{m\ep-1}+B_{v_{\min}}p\ep ,
\]
and the right-hand side tends to $0$ by \ref{A2} and \ref{A3}, so $v\ep\to\lambda/2$. If
$\lambda<0$, this contradicts the fact that $v\ep\ge v_{\min}>0$, which proves the last assertion of the
theorem. If $\lambda>0$, then $v\ep\ge\lambda/4$ for all sufficiently small $\varepsilon$, so
the argument above applies with $v_{\min}=\lambda/4$. Since $B_{\lambda/4}=B$, this gives the
stated bound. \hfill\qedsymbol
\begin{remark}
    The hypothesis $\liminf_{\epsilon\to0^+}v\ep>0$ is imposed to exclude bound states whose energies collapse to $0$ in the limit. Such states may occur for the regularized operators, but their energies approach the essential spectrum, which is $[0,\infty)$ for any Schro\"dinger operator with integrable potential~\cite[Lemma~9.35]{Teschl}, and in particular for every $H\ep$. 
\end{remark}

\subsection{Proof of Proposition~\ref{prop:B1}}
By Lemma~\ref{lem:conv}, $\psi\ep\in H^2(\mathbb{R})$, so by the Sobolev embedding theorem $\psi\ep\in C^1(\mathbb{R})$, with $\psi\ep$ and $\psi\ep\pri$ bounded and tending to $0$ as $\abs{x}\to\infty$.

(a) Multiplying \eqref{eq:bound} by $\overline{\psi\ep}$, integrating, and integrating the first term by parts (the boundary terms vanish since $\psi\ep,\psi\ep\pri\to0$ as $\abs{x}\to\infty$), we obtain, since $\norm{\psi\ep}_2=1$ by hypothesis,
\[
\norm{\psi\ep\pri}_2^2+v\ep^2=\lambda\int\delta\ep\abs{\psi\ep}^2\dd x,
\]
and hence, by \ref{A1},
\begin{equation}\label{norms}
\norm{\psi\ep\pri}_2^2+v\ep^2\leq\abs{\lambda}M\norm{\psi\ep}_\infty^2 .
\end{equation}
The lower bound claimed in (a) follows from \eqref{norms} whenever $v\ep\geq v_{\min}$. For the upper bound, Lemma~\ref{lem:elem}(d) gives
\begin{equation}\label{GNS}
\norm{\psi\ep}_\infty^2\leq2\norm{\psi\ep\pri}_2\norm{\psi\ep}_2=2\norm{\psi\ep\pri}_2 .
\end{equation}
Combining \eqref{norms} and \eqref{GNS} gives $\norm{\psi\ep\pri}_2^2\leq2\abs{\lambda}M\norm{\psi\ep\pri}_2$, so $\norm{\psi\ep\pri}_2\leq2\abs{\lambda}M$. Substituting this into \eqref{GNS} gives $\norm{\psi\ep}_\infty^2\leq4\abs{\lambda}M$, the upper bound in (a).

(b) Differentiating \eqref{conv} gives
\[
\psi\ep \pri (x)=\lambda\int G_{v\ep} \pri (x-y)\psi\ep(y)\delta\ep(y)\dd {y}.
\]
Since $\norm{G_{v\ep}\pri}_\infty=\frac12$ by \eqref{greens}, part (a) and \ref{A1} give
\[
\norm{\psi\ep \pri}_\infty \leq\abs{\lambda}\norm{G_{v\ep}\pri}_\infty\norm{\psi\ep}_\infty\norm{\delta\ep}_1
\leq \abs{\lambda}\cdot\frac12\cdot2\sqrt{\abs{\lambda}M}\cdot M
=(\abs{\lambda}M)^{3/2}.
\]
(c) By \eqref{norms} and the upper bound in (a), $v\ep^2\le\abs{\lambda}M\norm{\psi\ep}_\infty^2\le4\lambda^2M^2$.
\hspace*{\fill}\qedsymbol

\subsection{Proof of Proposition~\ref{prop:B2}}
(a) Using Lemma~\ref{lem:elem}(c) and Proposition~\ref{prop:B1},
\begin{align*}
    \frac{\abs{\lambda}}{2v\ep}\int \abs{e^{-v\ep\abs{y}}-1}\abs{\delta\ep(y)}\abs{\psi\ep(y)}\dd y
    &\leq\abs{\lambda}\sqrt{\abs{\lambda}M}\frac{(1+v\ep)}{v\ep}\int\abs{\delta\ep(y)}\min{\{1,\abs{y}\}}\dd y \\&\leq\abs{\lambda}\sqrt{\abs{\lambda}M}(1+\frac1{v_{\min}})p\ep.
\end{align*}
(b) By the fundamental theorem of calculus and Proposition~\ref{prop:B1}(b), 
\[
\abs{\psi\ep(y)-\psi\ep(0)}= \abs{\int_0^y\psi\ep\pri(x)\dd x}\leq\norm{\psi\ep\pri}_\infty\abs{y}\leq(\abs{\lambda}M)^{3/2}\abs{y}.
\]
Additionally, by the triangle inequality and Proposition~\ref{prop:B1}(a),
\[
\abs{\psi\ep(y)-\psi\ep(0)}\leq2\norm{\psi\ep}_\infty\leq4\sqrt{\abs{\lambda}M}.
\]
Hence, 
\[
    \abs{\psi\ep(y)-\psi\ep(0)}\leq\sqrt{\abs{\lambda}M}\min{\{4, \abs{\lambda}M\abs{y}\}}\leq4\sqrt{\abs{\lambda}M}(1+\abs{\lambda}M)\min{\{1, \abs{y}\}}.
\]
Thus,
\[
\frac{\abs{\lambda}}{2v\ep}\int \abs{\delta\ep(y)}\abs{\psi\ep(y)-\psi\ep(0)}\dd y\leq \frac2{v_{\min}}\abs{\lambda}\sqrt{\abs{\lambda}M}(1+\abs{\lambda}M)p\ep,
\]
as desired. \hfill \qedsymbol

\subsection{Proof of Proposition~\ref{prop:B3}}
Using the identity $\psi\ep(y)=\psi\ep(0)+\psi\ep(y)-\psi\ep(0)$ in \eqref{conv}, one has
\[
\psi\ep(x)=\frac\lambda{2v\ep}\int e^{-v\ep\abs{x-y}}\delta\ep(y)\psi\ep(0)\dd y + \frac\lambda{2v\ep}\int e^{-v\ep\abs{x-y}}\delta\ep(y)(\psi\ep(y)-\psi\ep(0))\dd y.
\]
Since $\abs{e^{-v\ep\abs{x-y}}}\leq1$, it follows from Proposition~\ref{prop:B2}(b) that the term containing $\psi\ep(y)-\psi\ep(0)$ is bounded in absolute value by $C_2p\ep$. Hence,
\[
\abs{\psi\ep(x)}\leq\frac{\abs{\lambda}M}{2v\ep}\abs{\psi\ep(0)}+C_2p\ep.
\]
Taking the supremum over $x$ and using Proposition~\ref{prop:B1}(a),
\[
\frac{v_{\min}}{\sqrt{\abs{\lambda}M}}\le\norm{\psi\ep}_\infty\le\frac{\abs{\lambda}M}{2v\ep}\abs{\psi\ep(0)}+C_2p\ep .
\]
Since $p\ep\to0$ by \ref{A3}, for all sufficiently small $\epsilon$ we have $C_2p\ep\le\frac{v_{\min}}{2\sqrt{\abs{\lambda}M}}$, and hence, using $v\ep\ge v_{\min}$,
\[
\abs{\psi\ep(0)}\ge\frac{2v\ep}{\abs{\lambda}M}\cdot\frac{v_{\min}}{2\sqrt{\abs{\lambda}M}}\ge\frac{v_{\min}^2}{(\abs{\lambda}M)^{3/2}}.\tag*{\qedsymbol} 
\]
\section{Conclusion}\label{sec:conclusion}
For the one-dimensional scattering problem with a regularized delta potential satisfying
\ref{A1}--\ref{A3}, we showed that the transmission and reflection coefficients converge to those of
the one-dimensional delta interaction, with rate $O(r\ep)$. This extends the known convergence of the
scattering data, previously established for families obtained by scaling a single fixed
profile~\cite[Thm.~I.3.4.1]{AGHH}, to a strictly larger class of regularizations. For bound states under \ref{A1}--\ref{A4}, we gave an
elementary and quantitative proof that any family of negative eigenvalues bounded away from $0$
converges to $-\lambda^2/4$ when $\lambda>0$, and showed that no such family exists when $\lambda<0$.
As explained in Section~\ref{ss:contribution}, this convergence can be obtained from~\cite[Thm.~3]{BFT};
our contribution to the bound-state problem is a new proof of this known result and the explicit rate
$O(r\ep)$ of convergence. Future work will seek to extend these results to point interactions in two
and three spatial dimensions, where regularizations of the delta interaction require a renormalization
to have a nontrivial limit.

\vspace{2mm}
\noindent\textbf{Author Contributions:} All authors contributed equally.

\smallskip
\noindent\textbf{Funding:} This research was supported in part by the Department of Mathematics and the College of Liberal Arts and Sciences, University of Florida.

\smallskip
\noindent\textbf{Data Availability Statement:} No new data were created or analyzed in this study.

\smallskip
\noindent\textbf{Conflicts of Interest:} The authors declare no conflicts of interest.
\begin{appendices}
\section{The Lippmann--Schwinger equation}\label{app:LS}
\begin{proof}[Proof of Lemma~\ref{lem:LS}]
Recall from Section~\ref{sec:scattering} that $\psi\ep\in C^1(\mathbb{R})$, $\psi\ep\pri\in AC_{\mathrm{loc}}$, and $\delta\ep\psi\ep\in L^1(\mathbb{R})$.

Using \eqref{eq:scat} we can rewrite the integral as
\begin{equation}\label{eq:integ}
\frac{i\lambda}{2k}\int e^{ik\abs{x-y}}\delta\ep(y)\psi\ep(y)\dd y=-\frac{i}{2k}\int e^{ik\abs{x-y}}(\psi\ep\dpri(y)+k^2\psi\ep(y))\dd y.
\end{equation}
Restricting the integration to $[-R,R]$, splitting at $y=x$, and integrating by parts twice, which is permissible since the kernel is smooth on each of $[-R,x]$ and $[x,R]$ and $\psi\ep,\psi\ep\pri\in AC_{\mathrm{loc}}$, we obtain
\begin{equation}\label{Rint}
\int^R_{-R} e^{ik\abs{x-y}}\big(\psi\ep\dpri(y)+k^2\psi\ep(y)\big)\dd y=e^{-ikx}\alpha_R-e^{ikx}\beta_R +2ik\psi\ep(x),
\end{equation}
where 
\[
\alpha_R=e^{ikR}[\psi\ep\pri(R)-ik\psi\ep(R)]\quad \text{and}\quad \beta_R=e^{ikR}[\psi\ep\pri(-R)+ik\psi\ep(-R)].
\]
To let $R\to\infty$, we first check that $\alpha_R$ and $\beta_R$ converge. Their derivatives with respect to $R$ are integrable:
\[
\frac{d}{dR}\alpha_R= e^{ikR}[\psi\ep\dpri(R)+k^2\psi\ep(R)]=-\lambda e^{ikR}\delta\ep(R)\psi\ep(R)\in L^1(\mathbb{R}),
\]
\[
\frac{d}{dR}\beta_R= e^{ikR}[-\psi\ep\dpri(-R)-k^2\psi\ep(-R)]=\lambda e^{ikR}\delta\ep(-R)\psi\ep(-R)\in L^1(\mathbb{R}).
\]
Hence, by the fundamental theorem of calculus, both limits exist and we can write
\[
\alpha \coloneqq\lim_{R\to\infty}{\alpha_R }\quad\text{and}\quad\beta\coloneqq\lim_{R\to\infty}\beta_R.
\]
By \eqref{eq:scat}, the left-hand side of \eqref{Rint} equals $-\lambda\int_{-R}^{R}e^{ik\abs{x-y}}\delta\ep(y)\psi\ep(y)\dd y$, which converges as $R\to\infty$ since $\delta\ep\psi\ep\in L^1(\mathbb{R})$. Taking $R\to\infty$ in \eqref{Rint}, substituting into \eqref{eq:integ} and rearranging, we obtain
\[
\psi\ep(x)=\frac{i\lambda}{2k}\int e^{ik\abs{x-y}}\delta\ep(y)\psi\ep(y)\dd y +\frac{i}{2k}\big[e^{-ikx}\alpha-e^{ikx}\beta\big].
\]
Letting $x\to+\infty$ and $x\to-\infty$ in this identity and matching with the asymptotic conditions in \eqref{eq:scat} forces $\alpha=0$ and $\beta=2ik$, which yields \eqref{eq:Lippmann} upon substitution. Moreover, by the definitions of $\alpha_R$ and $\beta_R$, $\alpha=0$ means that $\psi\ep\pri(x)-ik\psi\ep(x)\to0$ as $x\to+\infty$, and $\beta=2ik$ means that $\psi\ep\pri(x)+ik\psi\ep(x)=2ike^{ikx}+o(1)$ as $x\to-\infty$; combined with \eqref{eq:scat}, this gives \eqref{eq:psiprime}.
\end{proof}

\section{Weak convergence of the regularizations}\label{app:weak}
Assume \ref{A1}--\ref{A3}, let $f\in C_b(\mathbb{R})$, and fix $\eta\in(0,1)$. Splitting the integral at $|x|=\eta$ gives
\[
\abs{\int\delta\ep f(x)\dd x-f(0)}\leq \abs{f(0)}\abs{m\ep-1}+\int_{\abs{x}<\eta}\abs{\delta\ep}\abs{f(x)-f(0)}\dd x + \int_{\abs{x}\geq\eta}\abs{\delta\ep}\abs{f(x)-f(0)}\dd x.
\]
Since
$\min\{1,|x|\}\ge\eta$ for $|x|\ge\eta$, using \ref{A1}, we obtain
\[
\abs{\int\delta\ep f(x)\dd x-f(0)}\leq\norm{f}_\infty\abs{m\ep-1}+M\sup_{\abs{x}<\eta}{\abs{f(x)-f(0)}}+\frac{2\norm{f}_\infty}\eta p\ep.
\]
Letting $\epsilon\to0^+$, then $\eta\to0^+$, it follows that
\[
\int\delta\ep f(x)\dd x\too f(0). 
\]\hfill \qedsymbol

Note that since $\mathcal{D}(\mathbb{R})\subseteq C_b(\mathbb{R})$, this also shows that $\delta\ep\to\delta$ in $\mathcal{D}'(\mathbb{R})$. Thus, families satisfying \ref{A1}--\ref{A3} are in fact regularizations of the Dirac delta in the sense defined in the section \ref{sec:intro}. 
\end{appendices}

\end{document}